\documentstyle[11pt]{article}
\begin{document}
\newcommand{\p}{\parallel }
\makeatletter \makeatother
\newtheorem{th}{Theorem}[section]
\newtheorem{lem}{Lemma}[section]
\newtheorem{de}{Definition}[section]
\newtheorem{rem}{Remark}[section]
\newtheorem{cor}{Corollary}[section]
\renewcommand{\theequation}{\thesection.\arabic {equation}}

\title{{\bf Chern-Connes Character for the Invariant Dirac Operator in Odd Dimensions}
\thanks{partially supported by  MOEC and the 973 project.}}
\author{ Yong Wang \\
{\scriptsize \it Nankai Institute of Mathematics Tianjin 300071, P. R. China}\\
{\scriptsize \it ~~~~email: wangy581@nenu.edu.cn}}

\date{}
\maketitle

\begin{abstract}
 In this paper we give a proof of the
Lefschetz fixed point formula of Freed$^{\rm [1]}$ for an
orientation-reversing involution on an odd dimensional spin manifold
by using the direct geometric method introduced in [2] and then we
generalize this formula under the noncommutative geometry framework.
\\

\noindent{\bf 2000 MR Subject Classification}\quad 58j20
\end{abstract}

\section{ Introduction}

\quad In [2], Lafferty, Yu and Zhang presented a simple and direct
geometric proof of the Lefschetz fixed point formula for an
orientation-preserving isometry on an even dimensional spin
manifold by Clifford asymptotics of heat kernel. Chern and Hu [3]
used the method in [2] to compute the equivariant Chern-Connes
character for the invariant Dirac operator on an even dimensional
spin manifold. In [4] an alternate
approach to certain technical estimates in [3] was given.\\
 \indent In parallel, Freed$^{[1]}$ considered the case
of an orientation-reversing involution acting on an odd
dimensional spin manifold and gave the associated Lefschetz
formulas by K-theretical way$^{\rm [5]}$. The heat kernel
method$^{ [6]}$ may be used to prove this odd Lefschetz formula
as claimed in [1].\\
 \indent Inspired by [2] and [3], in this paper we
 give a direct geometric proof of the Freed's odd Lefschetz
formula. We also construct an even spectral triple (see Section 3)
by the Dirac operator and the orientation-reversing involution,
then compute
the Chern-Connes character for this spectral triple.\\
\indent The paper is organized as follows: In Section 2.1, we
present some notations and discuss the standard setup. Evaluating
the Clifford asymptotics of the local Lefschetz index is given in
Section 2.2 and Section 2.3. In Section 3, we construct an even
spectral triple and then compute its Chern-Connes character.

\section{ A direct geometric proof  of the Freed's odd Lefschetz
formula }

\noindent {\bf 2.1~ Preliminaries}\\

\indent Firstly we give the standard setup (also see Section 1 in
 [1]). Let $M$ be a closed, connected and oriented
Riemannian manifold of odd dimension $n$ with a fixed spin
structure ${\rm Spin}(M)$, and $S$ be the bundle of spinors on
$M$. Denote by $D$ the associated Dirac operator on $\Gamma(M;S)$,
the
 space of smooth sections of the bundle $S$. Let
$\tau:~M\rightarrow M$ be an orientation-reversing isometric
involution. Assume there exists a self-adjoint lift
$\widetilde{\tau}:~\Gamma(M;S)\rightarrow \Gamma(M;S)$ of $\tau$
satisfying
 $$\widetilde{\tau}^2=1;~~D\widetilde{\tau}=-\widetilde{\tau}D.\eqno(2.1)$$
When $\tau$ preserves Pin structure, such a lift
$\widetilde{\tau}$ always exists. Now the $+1$ and $-1$
eigenspaces of $\widetilde{\tau}$ give a splitting of the spinor
fields
$$\Gamma(M;S)\cong \Gamma^+(M;S){\small \oplus} \Gamma^{-}(M;S) \eqno(2.2)$$
and the Dirac operator interchanges $\Gamma^+(M;S)$ and
$\Gamma^-(M;S)$. We denote by $D^+$ the restriction of $D$ on
$\Gamma^+(M;S)$. The purpose of this section is to compute
$${\rm index}[D^+:~\Gamma^+(M;S)\rightarrow \Gamma^-(M;S)]\eqno(2.3)$$

In the following we give an explicit construction of
$\widetilde{\tau}$. The tangent map of $\tau$ gives a map
$d\tau:~O(M)\rightarrow O(M)$. Let the associated bundle ${\rm
Pin}(M)={\rm Spin}(M)\otimes_R{\rm Pin}(n)$ be the induced Pin
structure on $M$ where $R:~{\rm Spin}(n)\times {\rm
Pin}(n)\rightarrow {\rm Pin}(n)$ is the Clifford multiplication.
Assume $\tau$ preserves this Pin structure, i.e. $d\tau$ has a
lift $\overline{d\tau}$ such that the diagram
$$
\begin{array}{cc}
\ {\rm Pin}(M)   \\
 \   \pi\downarrow\\
 \  O(M)
\end{array}
\begin{array}{cc}
\ \begin{array}{cc} \  \overline{d\tau}\\
\ \longrightarrow \end{array}
  \\
 \  \\
 \ \longrightarrow\\
 \ \begin{array}{cc} \ {d\tau} \end{array}
\end{array}
\begin{array}{cc}
\ {\rm Pin}(M)   \\
 \ \downarrow\pi\\
 \  O(M)
\end{array}.
$$
is commutative where $\pi:{\rm Pin}(M)\rightarrow O(M)$ is the
double covering and $\overline{d\tau}$ commutes with the ${\rm
Pin}(n)$-action. We recall the odd dimensional Spin($n$)
representation$^{\rm [7]}$. Let ${\rm Cl}^+(n+1)$ be the even part
of the  Clifford algebra generated by $e_1,\cdots,e_{n+1}$ and {
$\cal {I_+}$} be the associated positive irreducible
representation. Let $\rho_1:{\rm Pin}(n)\rightarrow {\rm
Cl}^+(n+1);~e_i\rightarrow e_ie_{n+1}$ for $1\leq i\leq n$ be an
algebra homomorphism and $\rho_2:{\rm Cl}^+(n+1)\rightarrow {\rm
End}(\cal I_+)$ be a representation of ${\rm Cl}^+(n+1)^{[7]}$,
then $\rho=\rho_2\rho_1$ is a Pin($n$) representation. Note that
Spin$(n)$ is a subgroup of Pin$(n)$, so we have
$$S={\rm
Spin}(M)\times_{\rho}{\cal I_+}={\rm Pin}(M)\times_{\rho}{\cal
I_+}.\eqno(2.4)$$ A linear map $\widetilde{\tau_0}$ is defined as
follows. Suppose that $\phi\in\Gamma(S)$ is expressed locally over
an open set $U_{\tau x}$ by $\phi=[(\sigma,f)]$ for $x\in M$ and a
neighborhood $U_{\tau x}$ of $\tau x$, where $\sigma:~U_{\tau
x}\rightarrow {\rm Spin(M)}$ is a local Spin frame field and $f:~
U_{\tau x}\rightarrow {\cal {I_+}}$ is a spinor-valued function,
and $[(\sigma,f)]$ denotes the equivalence class of $(\sigma,f)$
in $S= {\rm Spin(M)}\times_{\rho}{\cal {I_+}}.$ Let
$$\widetilde{\tau_0}:~S_{\tau x}\rightarrow S_x;~(\widetilde{\tau_0}\phi)(x)
=[(((\overline{d\tau})^{-1}\sigma)(x),f(\tau x))] \eqno(2.5)$$
where the right of (2.5) denotes the equivalence class in $S={\rm
Pin}(M)\times_{\rho}{\cal {I_+}}$. Then
$\widetilde{\tau_0}D=-D\widetilde{\tau_0}$. By Lemma 1.5 of [1],
then $\widetilde{\tau_0}^2$ is a constant multiple of the
identity. Let $F_1,\cdots,F_r$ be components of the fixed point
set of $\tau$ and ${\rm codim}F_q=2m_q+1~ (1\leq q \leq r)$ and
$m_i\geq m_j$ for $i<j$, then $\widetilde{\tau_0}^2=(-1)^{m_1+1}$
over the neighborhood of $F_1$ (see Section 2.2). So we define
$$\widetilde{\tau}=(\sqrt{-1})^{m_1+1}\widetilde{\tau_0},\eqno(2.6)$$
then $\widetilde{\tau}$ satisfies the condition (2.1). Note that
since $\tau$ preserves the Pin structure, ${\rm
codim}F_i\equiv{\rm codim}F_j~{\rm mod}~ 4$ (similar to
  Proposition 8.46 in [8]). So (2.6) up to a
sign is independent of the choice of components. We take the
Pin($n$)-invariant Hermitian inner product on ${\cal {I_+}}$, then
by (2.5) and (2.6), we have
$\widetilde{\tau}{\widetilde{\tau}}_a=1$ where
${\widetilde{\tau}}_a$ is the adjoint operator of
$\widetilde{\tau}$. Considering $\widetilde{\tau}^2=1$ then
$\widetilde{\tau}={\widetilde{\tau}}_a$.\\ \indent By
Mckean-Singer formula, we have
$${\rm Ind}D^+={\rm Tr}(\widetilde{\tau}e^{-tD^2}).\eqno(2.7)$$
Let $P_t(x,y):S_y\rightarrow S_x$ be the fundamental solutions for
the heat operator $\partial/\partial t+D^2.$ The standard heat
equation argument yields
$${\rm Tr}(\widetilde{\tau}e^{-tD^2})=\int_M{\rm
Tr}[\widetilde{\tau}P_t(\tau x,x)]dx. \eqno(2.8)$$ We shall use
the abbreviation ${\cal{L}}(t,x)={\rm Tr}[\widetilde{\tau}P_t(\tau
x,x)]$. Let $\nu$ be the normal bundle of the fixed point set and
${\nu}(\varepsilon)=\{ x\in{\nu}|~||x||<\varepsilon\}$ for
$\varepsilon>0$. Similar to the discussions in [2], we get

\indent {\bf Theorem 2.1}
$${\rm Ind}D^+=\sum_{q=1}^r\int_{F_q}{\cal{L}}_{\rm loc}(\tau)(\xi)d\xi \eqno(2.9)$$
{\it where}$$ {\cal{L}}_{\rm loc}(\tau)(\xi)={\rm
lim}_{t\rightarrow 0}\int_{\nu_{\xi}(\varepsilon)}{\cal{L}}(t,{\rm
exp}c)dc\eqno(2.10)$$ {\it exists and is independent of}
$\varepsilon$.

 Since $\tau$ preserves the Pin structure, each $F_q$ has a natural orientation
 (similar to
Proposition 6.14 in [6]). Let ${\rm dim}F_q=2n'$ and $\xi\in F_q$,
then as in [2] there exists an oriented orthonormal frame field
$E=(E_1,\cdots,E_n)$
in a neighborhood $U$ of $\xi$ such that:\\
~(a)~ for $\zeta\in U\bigcap F_q$,
$(E_1(\zeta),\cdots,E_{2n'}(\zeta))$ is an oriented orthonormal
basis of $T_{\zeta}F_q$ while the vector
fields $E_{2n'+1}(\zeta),\cdots,E_{n}(\zeta)$ are normal to $T_{\zeta}F_q$.\\
 ~(b)~ $E$ is parallel along the geodesics normal to $F_q$.\\
With respect to $(E_1,\cdots,E_n)$, $d\tau$ is expressed as a
matrix-valued function ${\cal T}$ for $x\in U$
$$d\tau E(x)=E(\tau x){\cal T}(x).$$
Moreover there is a neighborhood $V$ of $\xi$ in $F_q$ such that
$E$ is defined on $U={\rm exp}(\nu|_V\bigcap \nu(\varepsilon))$
for sufficient small $\varepsilon$. If $B_0(\varepsilon)$ is the
ball of radius $\varepsilon$ in ${\bf R}^{2m_q+1}$; we define the
homeomorphism $\Phi:~V\times B_0(\varepsilon)\rightarrow U$ by
setting
$$\Phi(x';c_1,\cdots,c_{2m_q+1})=x={\rm exp}_{x'}(\sum^{2m_q+1}_{\alpha=1}c_{\alpha}E_{2n'+\alpha}(x')).\eqno(2.11)$$
Denote by $(x';c)$ the orthogonal coordinates of $x$ with respect
to $E=(E_1,\cdots,E_n)$ at $\xi$. Then we have that\\
~(a)  ${\cal T}(x';c)={\cal T}(x')$; and\\
~(b)  the isometry $\tau$ has the form $\tau(x';c)=(x';-c)$ and
for $\forall x\in U$
$${\cal T}(x)=\left[\begin{array}{lcr}
  \ I & 0  \\
    \  0  & -I
\end{array}\right];
$$
~(c)  Let $E^{\tau x}$ be an oriented frame field defined over the
patch $U$ by requiring that $E^{\tau x}(\tau x)=E(\tau x)$ and
that $E^{\tau x}$ be parallel along geodesic through $\tau x$.
Define the coordinates $\{ y_i\}$ of $x$ as
$$ (x';c)=x={\rm exp}_{\tau x}(\sum^n_{i=1}y_iE^{\tau x}_i(\tau
x))$$ then $E^{\tau x}(x)=E(x)$ and $y_i=0$ for $1\leq i \leq
2n';$ $y_{2n'+\alpha}=2c_\alpha$ for $1\leq \alpha \leq 2m_q+1.$
Note that (b) comes from $\tau^2=$id, i.e. ${\cal T}(x)^2=$id.
Since $x=(x',c),~x'$ and $(x',-c)=\tau x$ belong to the same
geodesic normal to $F_q$, (c) is correct.

\noindent {\bf 2.2 The Clifford asymptotics}

 Choose a Spin frame field $\sigma:U\rightarrow {\rm Spin}(M)$ such
 that $\pi'\sigma=(E^{\tau x}_1,\cdots,E^{\tau x}_n)$ where $\pi':
 {\rm Spin}(M)\rightarrow SO(M)$ is the double covering. For $x\in U$,
 let $\overline{P_t}(x)$, $\widetilde{\tau}^*(x)\in {\rm
 Hom}({\cal {I_+}},{\cal {I_+}})$ be defined through the equivalence relations for
 $u,v\in {\cal {I_+}}$
 $$ P_t(\tau x,x)[(\sigma(x),v)]=[(\sigma(\tau x),\overline{P_t}(x)v)]\eqno(2.12)$$
and
$$\widetilde{\tau}[(\sigma(\tau x),u)]=[(\sigma(x),\widetilde{\tau}^*(x)u].\eqno(2.13)$$
Similar to Lemma 4.1 in [10], we have:

 \indent {\bf Lemma
2.2}~~{\it For $x$ in a sufficient small neighborhood of $F_q$ and
$t>0$, the integrand ${\cal L}(t,x)$ is evaluated by}
$${\cal L}(t,x)={\rm
Tr}(\widetilde{\tau}^*(x)\overline{P_t}(x)).\eqno(2.14)$$ \indent
As in [2] and [9], in the normal coordinates $y_1,\cdots,y_n$ at
$\tau x$ with respect to the frame field $E^{\tau x}=(E^{\tau
x}_1,\cdots,E^{\tau x}_n)$, the operator
$${\chi}(y^{\alpha}D^{\beta}_ye^{\gamma})=|\beta|-|\alpha|+|\gamma|\eqno(2.15)$$
for multi-indices $\alpha,~\beta,$ and $\gamma$, with
$y^{\alpha}=y^{\alpha_1}_1\cdots y^{\alpha_n}_n,$
$D^{\beta}_y=(\partial/\partial_1)^{\beta_1}\cdots
(\partial/\partial_n)^{\beta_n}$ and
$e^{\gamma}=e^{\gamma_1}_1\cdots e^{\gamma_n}_n$ for
$\gamma_i\in\{ 0,1 \}.$

In the coordinates $(x_0;c)$ with respect to the frame field
$E=(E_1,\cdots,E_n)$, set $c=\sqrt{t}b$ and define operator
$\overline{\chi}$ on the monomials $\phi(t)e_{i_1}\cdots e_{i_s}$
by $$\overline{\chi}(\phi(t)e_{i_1}\cdots e_{i_s})=s-{\rm sup}\{
l\in {\bf Z}{\bf |}~{\rm lim}_{t\rightarrow
0}\frac{|\phi(t)|}{t^{l/2}}<\infty \},\eqno(2.16)$$ where
$\phi(t)\in {\bf R}.$ We denote by $P=Q+(\overline{\chi}<m)$ the
congruence of $P$ and $Q$ modulo the space generated by monomials
with $\overline{\chi}<m$.

By Section 2.1, then
$$d\tau E^{\tau x}(x)=d\tau E(x)=E(\tau x){\cal T}(x)=E^{\tau x}(\tau x)\left[\begin{array}{lcr}
  \ I & 0  \\
   \  0  & -I
\end{array}\right].\eqno(2.17)$$ Thus
$$\overline{d\tau}\sigma(x)=\sigma(\tau
x)c(e_{2n'+1})\cdots c(e_n).\eqno(2.18).$$ By (2.6), then
$$\widetilde{\tau}^*=(\sqrt{-1})^{m_1+1}c(e_{2n'+1})\cdots
c(e_n)=(\sqrt{-1})^{m_1-m_q}((\sqrt{-1})^{m_q+1}c(e_{2n'+1})\cdots
c(e_n)).\eqno(2.19)$$ \indent Let $x'$ be a point near $x=(\xi;c)$
and let
 $y=(y'_1,\cdots,y'_n)$ be the normal coordinate of $x'$ at point
 $\tau x$ with respect to the orthonormal frame field $E^{\tau x}=(E^{\tau x}_1,\cdots,E^{\tau
 x}_n)$ defined in Section 2.1. As in [2], let $\widetilde{A}$ be
 the $n\times n$ matrix defined by
 $$\widetilde{A_{ij}}=-\frac{1}{2}\sum^n_{k,l=1}R^{\tau x}_{ijkl}(\tau
 x)c(e_k)c(e_l),\eqno(2.20)$$
where $R^{\tau x}_{ijkl}(\tau
 x)$ are the coefficients of the Riemannian curvature tensor under
 the frame field $E^{\tau x}=(E^{\tau x}_1,\cdots,E^{\tau
 x}_n)$ at point $\tau x$. We define $\widetilde{A}^l(y')$ as
$${\widetilde{A}}^l(y')=\sum^n_{i,j=1}y'_iy'_j{\widetilde{A}}^l_{ij}\eqno(2.21)$$
for $l=1,2,\cdots.$

Similarly to [9], in the odd dimensional case, there is a function
$P(t;z_1,z_2,\cdots;w_1,w_2,\cdots)$, which is a power series in
$t$ with coefficient polynomials in $z_i$ and
$w_i$ such that\\
 $~~\overline{P_t(x)}=(1/4\pi t)^{\frac{n}{2}}{\rm exp}(-d^2(x,\tau
 x)/4t)$
$$\times[P(t;{\rm Tr}\widetilde{A}^2,\cdots,{\rm Tr}\widetilde{A}^{2k},\cdots;
\widetilde{A}^2(y)\cdots,\widetilde{A}^{2l}(y),\cdots)+\sum_{m\geq
0 }t^m(\overline{\chi}<2m)],\eqno(2.22)$$ where in the diagonal
form we have, by solving harmonic oscillator-type equations,
$$P(t;((-1)^k2(x^{2k}_1+\cdots+x^{2k}_{n'+m_q}));((-1)^l\sum_{\alpha=1}^{n'+m_q}(y_{2\alpha-1}^2
+y_{2\alpha}^2)x_{\alpha}^{2l}))=(4\pi t)^{n'+m_q}$$
$$\times {\rm
exp}(\sum^{n-1}_{\alpha=1}y_{\alpha}^2/4t)\prod_{s=1}^{n'+m_q}\left[
\frac{\sqrt{-1}x_s}{8\pi {\rm sinh}\frac{\sqrt{-1}x_st}{2}}
 {\rm exp}(-\frac{\sqrt{-1}x_s}{8}{\rm coth}\frac{\sqrt{-1}x_st}{2}(y_{2s-1}^2+y_{2s}^2))
 \right].\eqno(2.23)$$
\indent As Lemma 4.3 in [2], we have:

 \indent {\bf Lemma 2.3}~~{\it Let
$\widetilde{A_0}$ be
 the matrix defined by
 $$\widetilde{A_0}_{ij}=-\frac{1}{2}\sum^{2n'}_{k,l=1}R_{ijkl}(\xi)c(e_k)c(e_l),~~~1\leq i,j\leq n ,\eqno(2.24)$$
and define the tangential component $A^{\top}$ and the normal
component $A^{\bot}$ by
 $${A^{\top}}_{ij}=-\frac{1}{2}\sum^{2n'}_{k,l=1}R_{ijkl}(\xi)c(e_k)c(e_l),~~~1\leq i,j\leq 2n', \eqno(2.25)$$
$${A^{\bot}}_{ij}=-\frac{1}{2}\sum^{2n'}_{k,l=1}R_{ijkl}(\xi)c(e_k)c(e_l),~~~2n'+1\leq i,j\leq n .\eqno(2.26)$$
\noindent Then
$$\widetilde{A_0}=\left(\begin{array}{lcr}
  \ A^{\top} & 0  \\
      ~~0  & A^{\bot}
\end{array}\right).\eqno(2.27)$$
Further, the relations
$${\rm Tr}\widetilde{A}^{2k}={\rm Tr}(A^{\top})^{2k}+{\rm
Tr}(A^{\bot})^{2k}+\sum_{\alpha=1}^{n-2n'}c(e_{2n'+\alpha})(\overline{\chi}<4k)+(\overline{\chi}<4k)\eqno(2.28)$$
and
$$\widetilde{A}^{2k}(y/\sqrt{t})=4(A^{\bot})^{2k}(b)+
\sum_{\alpha=1}^{n-2n'}c(e_{2n'+\alpha})(\overline{\chi}<4k)+(\overline{\chi}<4k)\eqno(2.29)$$
hold.}

Combining (2.19),(2.22) and Lemma 2.3, we get

\indent {\bf Lemma 2.4}
$$\widetilde{\tau}^*(x)\overline{P_t}(x)=e^{-||b||^2}\left[\frac{(\sqrt{-1})^{m_1+1}}{(4\pi t)^{n/2}}
P(t;({\rm Tr}(A^{\top})^{2k}+{\rm
Tr}(A^{\bot})^{2k});(4t(A^{\bot})^{2k}(b)))\right.$$
$$\left.\times c(e_{2n'+1})\cdots
c(e_n)+(\overline{\chi}<2n-2n')_b\right]\eqno(2.30)$$ {\it where
$c=\sqrt{t}b$ and $(\overline{\chi}<2n-2n')_b$ denotes the space
spanned by which are polynomials in $b$ and satisfy
$\overline{\chi}<2n-2n'$.}

\noindent {\bf 2.3 Evaluation of the local index}

\indent {\bf Lemma 2.5}~([2])
$${\rm lim}_{t\rightarrow 0}\int _{\nu_{\xi}(\varepsilon)}e^{-||b||^2}{\rm
Tr}(\phi) dc=0\eqno(2.31)$$ {\it where} $\phi\in (\overline{\chi}<2n-2n')_b.$\\
\indent To compute the trace it suffices to compute the
coefficient of the $c(e_1)\cdots c(e_n)$ term in Lemma 2.4. Note
that $A^{\bot}$ and $A^{\top}$ are of order $\overline{\chi}\leq
2$, containing terms $c(e_i)c(e_j)$ with $1\leq i,j \leq 2n'$ and
${\rm Tr}(c(e_1)\cdots
c(e_n))=(-\sqrt{-1})^{[\frac{n}{2}]+1}2^{[\frac{n}{2}]}$. Since
$c(e_i)c(e_j)=-c(e_j)c(e_i)+(\overline{\chi}<1)$, if we formally
replace $c(e_i)$ by $\omega_i$ where
$\omega=(\omega_1,\cdots,\omega_n)$ is the frame dual to $E$, and
then substitute $\Omega^{\top}$ and $\Omega^{\bot}$ for $A^{\top}$
and $A^{\bot}$, where
$$\Omega^{\top}=-\frac{1}{2}\sum^{n}_{k,l=1}R_{ijkl}(\xi)\omega_k\wedge\omega_l,~~~1\leq i,j\leq 2n' ;\eqno(2.32)$$
$${\Omega^{\bot}}=-\frac{1}{2}\sum^{n}_{k,l=1}R_{ijkl}(\xi)\omega_k\wedge\omega_l,~~~2n'+1\leq i,j\leq n. \eqno(2.33)$$
To compute the trace, we only need to compute the top form (of
order $2n'$) on $F_q$, then we multiply it by
$(-\sqrt{-1})^{[\frac{n}{2}]+1}2^{[\frac{n}{2}]}$. In order to
compute this differential form, we need the odd dimensional case
of the Chern root algorithm (see [9]).

Let $$\Omega=\left[\begin{array}{lcr}
  \ \Omega^{\top} & 0  \\
      ~~0  & \Omega^{\bot}
\end{array}\right]$$
be given formally as $$\Omega^{\top}=\left[\begin{array}{lllll}
0 & u_1 & & & \\
-u_1 & 0 & & & \\
& &\ddots & & \\
& & & 0 & u_{n'} \\
& & & -u_{n'} & 0 \\
\end{array}\right],~~~~~~
\Omega^{\bot}=\left[\begin{array}{llllll}
0 & v_1 & & & & \\
-v_1 & 0 & & & & \\
& &\ddots & & & \\
& & & 0 & v_{m_q} & \\
& & & -v_{m_q} & 0 & \\
& & & & & 0 \\
\end{array}\right],$$
where $u_i$ and $v_i$ are indeterminants. Then
$$4t(\Omega^{\bot})^{2k}(b)=(-1)^k4t\sum_{\alpha=1}^{m_q}v_{\alpha}^{2k}(b^2_{2\alpha-1}+b^2_{2\alpha});\eqno(2.34)$$
$${\rm Tr}\Omega^{2k}=2(-1)^k(\sum_{\alpha=1}^{n'}u^{2k}_{\alpha}+\sum_{\beta=1}^{m_q}v_{\beta}^{2k}).\eqno(2.35)$$
By (2.23),(2.30),(2.34) and (2.35), we have:
\begin{eqnarray*}
{\cal L}_{\rm loc}(\tau)&=&{\rm lim}_{t\rightarrow 0}\int _{{\bf R
}^{n-2n'}}(-\sqrt{-1})^{[\frac{n}{2}]+1}2^{[\frac{n}{2}]}t^{\frac{n-2n'}{2}}e^{-||b||^2}
\frac{(\sqrt{-1})^{m_1+1}}{(4\pi t)^{n/2}}\\
& ~ &\times
P(t;(2(-1)^k(\sum_{\alpha=1}^{n'}u^{2k}_{\alpha}+\sum_{\beta=1}^{m_q}v_{\beta}^{2k}));
((-1)^k4t\sum_{\alpha=1}^{m_q}v_{\alpha}^{2k}(b^2_{2\alpha-1}+b^2_{2\alpha})))db\\
& = & {\rm lim}_{t\rightarrow 0}\int _{{\bf R
}^{n-2n'}}(-\sqrt{-1})^{[\frac{n}{2}]+1}2^{[\frac{n}{2}]}t^{\frac{n-2n'}{2}}e^{-b_n^2}
\frac{(\sqrt{-1})^{m_1+1}}{(4\pi t)^{n/2}}\\
& ~ &
\times\prod^{n'}_{\alpha=1}\frac{\sqrt{-1}tu_{\alpha}/2}{{\rm
sinh}\sqrt{-1}tu_{\alpha}/2}
\prod^{m_q}_{\beta=1}\frac{\sqrt{-1}tv_{\beta}/2}{{\rm
sinh}\sqrt{-1}tv_{\beta}/2}\\
& ~ & \times {\rm exp}(-\sum
^{m_q}_{s=1}\frac{\sqrt{-1}v_st}{2}{\rm
coth}\frac{\sqrt{-1}v_st}{2} (b_{2s-1}^2+b_{2s}^2))db.
\end{eqnarray*}\\
Note that $\int_{\bf R}e^{-b_n^2}db_n=\sqrt{\pi}$. In the final
calculation after integrating out $b$, we will take the form of
order $2n'$ on $F_q$, and hence the factor of $t$ cancels. So
\begin{eqnarray*}
{\cal L}_{\rm loc}(\tau)&=&\sqrt{\pi} \int _{{\bf R
}^{2m_q}}(-\sqrt{-1})^{[\frac{n}{2}]+1}2^{[\frac{n}{2}]}
\frac{(\sqrt{-1})^{m_1+1}}{(4\pi)^{n/2}}\\
&~&\times\prod^{n'}_{\alpha=1}\frac{\sqrt{-1}u_{\alpha}/2}{{\rm
sinh}\sqrt{-1}u_{\alpha}/2}
\prod^{m_q}_{\beta=1}\frac{\sqrt{-1}v_{\beta}/2}{{\rm
sinh}\sqrt{-1}v_{\beta}/2}\\
& ~ & \times {\rm exp}(-\sum
^{m_q}_{s=1}\frac{\sqrt{-1}v_s}{2}{\rm coth}\frac{\sqrt{-1}v_s}{2}
(b_{2s-1}^2+b_{2s}^2))db_1\cdots db_{2m_q}\\
&=&(-\sqrt{-1})^{[\frac{n}{2}]+1}2^{[\frac{n}{2}]}\sqrt{\pi}
\frac{(\sqrt{-1})^{m_1+1}}{(4\pi)^{n/2}}\int _{{\bf R }^{2m_q}}
\prod^{n'}_{\alpha=1}\frac{\sqrt{-1}u_{\alpha}/2}{{\rm
sinh}\sqrt{-1}u_{\alpha}/2}
\prod^{m_q}_{\beta=1}\frac{\sqrt{-1}v_{\beta}/2}{{\rm
sinh}\sqrt{-1}v_{\beta}/2}\\
&~&\times {\rm exp}(-\frac{1}{2}\sum_s\frac{v_s}{{\rm
sin}v_s/2}{\rm sin}(\frac{\pi+v_s}{2})(b_{2s-1}^2+b_{2s}^2))db_1\cdots db_{2m_q}\\
&=&\frac{(-\sqrt{-1})^{[\frac{n}{2}]+1}(\sqrt{-1})^{m_1+1}}
{\pi^{n'}2^{[\frac{n}{2}]+1}}
\prod^{n'}_{\alpha=1}\frac{\sqrt{-1}u_{\alpha}/2}{{\rm
sinh}\sqrt{-1}u_{\alpha}/2} \prod^{m_q}_{\beta=1}[{\rm
sin}(\frac{v_\beta+\pi}{2})]^{-1}.
\end{eqnarray*}
Let
$\frac{u_{\alpha}}{2\pi}=u_{\alpha}^*,~\frac{v_{\beta}}{2\pi}=v_{\beta}^*$
be the Chern roots, then
 \begin{eqnarray*}
{\cal L}_{\rm loc}(\tau)^{(2n')}&=&\left[
\frac{(-\sqrt{-1})^{[\frac{n}{2}]+1}(\sqrt{-1})^{m_1+1}}{\pi^{n'}2^{[\frac{n}{2}]+1}}(\sqrt{-1})^{n'}(2\pi)^{n'}\right.\\
& &\times\left.\prod^{n'}_{\alpha=1}\frac{u_{\alpha}/4\pi}{{\rm
sinh}u_{\alpha}/4\pi} \prod^{m_q}_{\beta=1}[{\rm
sin}(\frac{v_\beta}{4\pi\sqrt{-1}}+\frac{\pi}{2})]^{-1}\right]^{(2n')}\\
&=&\left[\frac{(\sqrt{-1})^{m_1}}{2}
\prod^{n'}_{\alpha=1}\frac{u_{\alpha}^*/2}{{\rm
sinh}u_{\alpha}^*/2} \prod^{m_q}_{\beta=1}[2{\rm
sinh}(\frac{v_\beta}{4\pi}+\frac{\sqrt{-1}\pi}{2})]^{-1}\right]^{(2n')}\\
&=&\frac{(\sqrt{-1})^{m_1-m_q}}{2}\left[
\prod^{n'}_{\alpha=1}\frac{u_{\alpha}^*/2}{{\rm
sinh}u_{\alpha}^*/2}
\prod^{m_q}_{\beta=1}(e^{\frac{v^*_{\beta}}{2}}+e^{-\frac{v^*_{\beta}}{2}})^{-1}\right]^{(2n')}.
\end{eqnarray*}
As in [5], we write the characteristic class
$${\rm ch}\triangle(N_q)=\prod^{m_q}_{\beta=1}(e^{\frac{v^*_{\beta}}{2}}+e^{-\frac{v^*_{\beta}}{2}})\eqno(2.36)$$
where $N_q$ denotes the normal bundle of $F_q$. We thus obtain the following theorem. \\

\indent {\bf Theorem 2.6}~([1])~{\it Let $M$ be an odd dimensional
compact oriented Spin manifold and $\tau:M\rightarrow M$ be an
orientation-reversing isometric involution which preserves Pin
structure. Suppose that $F_1,\cdots,F_r$ are components of the
fixed point set, then}
$${\rm ind}D^+=\frac{1}{2}\sum_{q=1}^r
\int_{F_q}(\sqrt{-1})^{m_1-m_q}\widehat{A}(TF_q)[{\rm
ch}\triangle(N_q)]^{-1}.\eqno(2.37)$$

\indent {\bf Remark}~By (2.6), the grading operator
$\widetilde{\tau}$ depends on $m_1$. If we choose another
component $F_i$, then $\widetilde{\tau}$ is up to
$(\sqrt{-1})^{m_i-m_1}$, but we note by the change of
$\widetilde{\tau}$ and (2.2),(2.3), then ${\rm ind}D^+$ also
change $(\sqrt{-1})^{m_i-m_1}.$

\section{ The Chern-Connes character of even spectral triple
$(C^{\infty}_{\tau}(M),L^2(M,S),D,\widetilde{\tau})$}

\quad Let $M,\tau$ and $\widetilde{\tau}$ be given as in Section
2.1, let
$$C^{\infty}_{\tau}(M)=\{ a\in
C^{\infty}(M)|a\tau(x)=a(x);~\forall x\in M \}\subset
C^{\infty}(M),$$ then
$(C^{\infty}_{\tau}(M),L^2(M,S),D,\widetilde{\tau})$ is an
$\theta$-summable even spectral triple (for definition see [10] or
[11] ). In the following we will compute its Chern-Connes
character. Firstly let us review the definition of the
Chern-Connes character represented by the JLO cocycle
 in the entire cyclic cohomology .

\indent {\bf Definition 3.1}~([12])~Let $(A,H,D,\gamma)$ be an
even $\theta$-summable spectral triple associated to a Banach
algebra $A$ with identity, then its Chern character ${\rm
ch}_{*}(A,H,D,\gamma)=\{ {\rm ch}_k(D)|~k\geq 0 ~{\rm and~ even}
\}$ in the entire cyclic cohomology is defined by
$${\bf {\rm
ch}}_k(D)(a^0,\cdots ,a^k)=\int_{\triangle_k}{\rm str}
(a^0e^{-s_1D^2}[D,a^1]e^{-(s_2-s_1)D^2}[D,a^2]$$ $$\cdots
e^{-(s_k-s_{k-1})D^2}[D,a^k] e^{-(1-s_k)D^2})ds,\eqno(3.1)$$ where
$a^i\in A$ and $\triangle_k=\{(s_1,\cdots,s_k)|~0\leq
s_1\leq\cdots\leq s_k\leq 1\}$. For $t>0$, considering the
deformed Chern-Connes character ${\rm ch}_{*}(\sqrt{t}D)=\{ {\rm
ch}_k(\sqrt{t}D)|~k\geq 0 ~{\rm and~ even}\}$ is expressed by
$${\bf {\rm
ch}}_k(\sqrt{t}D)(a^0,\cdots
,a^k)=t^{\frac{k}{2}}\int_{\triangle_k}{\rm str}
(a^0e^{-s_1tD^2}[D,a^1]e^{-(s_2-s_1)tD^2}[D,a^2]$$ $$\cdots
e^{-(s_k-s_{k-1})tD^2}[D,a^k] e^{-(1-s_k)tD^2})ds.\eqno(3.2)$$
 We write
$$\lambda(p)=(\lambda_1,\cdots,\lambda_p);~|\lambda(p)|=\lambda_1+\cdots+\lambda_p;
~\lambda(p)!=\lambda_1!\cdots\lambda_p!$$
$$\widetilde{\lambda}(p)!=(\lambda_1+1)(\lambda_1+\lambda_2+2)\cdots(\lambda_1+\cdots+\lambda_{p}+p).$$
For an operator $B$ and any positive integer $l$, write
$B^{[l]}=[D^2,B^{[l-1]}],~B^{[0]}=B$. We use the notation
$$D^{\lambda(p)}=f^0[c(df^{1})]^{[\lambda _{1}]}
[c(df^2)]^{[\lambda _2]}\cdots [c(df^p)]^{[\lambda
_p]};~~D^{\lambda(p)}_t=t^{\frac{p}{2}+|\lambda(p)|}D^{\lambda(p)}\eqno(3.3)$$
where $f^j\in C^{\infty}_{\tau}(M)$ for $0\leq j\leq p$. By $0\leq
\lambda(p)\leq n-k$, we mean $0\leq \lambda_j\leq n-k$ for $1\leq
j\leq p$. Recall a result in [3] or [4].

 \indent {\bf Lemma
3.2}~([3],[4])~(i)~{\it When $k\leq n$ and $t\rightarrow 0^+$, we
have:
$${\bf {\rm ch}}_k(\sqrt
{t}D)(f^0,\cdots ,f^k)=\sum_{0\leq \lambda(k)\leq
n-k}\frac{(-1)^{|\lambda(k)|}t^{|\lambda(k)|+\frac{k}{2}}}{\lambda(k)!\widetilde{\lambda(k)}!}
{\rm str}\{D^{\lambda(k)}e^{-tD^2}\}+O(\sqrt{t})¡£\eqno(3.4)$$}
~(ii)~{\it If $k>n$, then when $t\rightarrow 0^+$, we have:}
$${\rm lim}_{t\rightarrow 0}{\bf {\rm ch}}_k(\sqrt {t}D)(f^0,\cdots
,f^k)=0.\eqno(3.5)$$
 \indent In the following, we'll compute ${\rm lim}_{t\rightarrow 0}{\rm
 str}\{D^{\lambda(k)}_te^{-tD^2}\}$ by using the method in
 Section 2. We consider the coordinates systems in Section 2.1,
 2.2. Similar to Theorem 2.1 and Lemma 2.2, we have

\indent {\bf Proposition 3.3}~ $${\rm lim}_{t\rightarrow 0}{\rm
 Tr}(\widetilde{\tau}D^{\lambda(k)}_te^{-tD^2})=\sum^r_{q=1}\int_{F_q}
 ({\rm
lim}_{t\rightarrow 0}\int_{\nu_{\xi}(\varepsilon)}{\rm
Tr}(\widetilde{\tau}^*(x)\overline{P_t^{\lambda}}(x))dx)d\xi,\eqno(3.6)$$
{\it where}
$$D^{\lambda(k)}_tP_t(\tau x,x)[(\sigma(x),v)]=[(\sigma(\tau x),\overline{P_t^{\lambda}}(x)v)].\eqno(3.7)$$

\indent {\bf Lemma 3.4}~([3])~ {\it For $f^0,\cdots,f^k\in
C^{\infty}_{\tau}(M)$ and $\lambda\neq 0$, then}
$\overline{\chi}(D^{\lambda(k)}_t)<0$.

Let $\xi\in F_q$, we will compute ${\rm lim}_{t\rightarrow
0}\int_{\nu_{\xi}(\varepsilon)}{\rm
Tr}(\widetilde{\tau}^*(x)\overline{P_t^{\lambda}}(x))dx.$

\indent {\bf Theorem 3.5}~~{\it If
$\lambda=(\lambda_1,\cdots,\lambda_k)\neq 0$, then we have:}
$$\widetilde{\tau}^*(x)\overline{P_t^{\lambda}}(x)=e^{-||b||^2}(\overline{\chi}<2n-2n')_b,\eqno(3.8)$$
{\it and} ${\rm lim}_{t\rightarrow
0}\int_{\nu_{\xi}(\varepsilon)}{\rm
Tr}(\widetilde{\tau}^*(x)\overline{P_t^{\lambda}}(x))dx=0$ {\it
i.e.} ${\rm lim}_{t\rightarrow 0}{\rm
 Tr}(\widetilde{\tau}D^{\lambda(k)}_te^{-tD^2})=0$.
 \\
 {\it Proof.}~~This theorem comes from Lemma 2.4, Lemma 3.4
and Lemma 2.5.~~~$\Box$\\
\indent As in [3], for any $g\in
C^{\infty}_{\tau}(M),~x=(\xi;c)=(\xi;\sqrt{t}b)$, then
$g(x)=g(\xi)+(\overline{\chi}<0)$ and
$(dg)(x)=(dg)(\xi)+(\overline{\chi}<1).$ So we have
$$t^{\frac{k}{2}}f^0c(df^1)\cdots c(df^k)=t^{\frac{k}{2}}f^0(\xi)c(df^1)(\xi)\cdots
c(df^k)(\xi)+(\overline{\chi}<0).\eqno(3.9)$$ By Lemma 2.4 and
(3.9), we obtain:
$$\widetilde{\tau}^*(x)t^{\frac{k}{2}}f^0c(df^1)\cdots c(df^k)\overline{P_t}(x)
=(\sqrt{-1})^{m_1+1}t^{\frac{k}{2}}f^0(\xi)c(d(f^1))(\xi)\cdots
c(d(f^k))(\xi)$$ $$\times e^{-||b||^2}\frac{1}{(4\pi t)^{n/2}}
P(t;(\cdots,({\rm Tr}(A^{\top})^{2l}+{\rm
Tr}(A^{\bot})^{2l}),\cdots);(\cdots,(4t(A^{\bot})^{2l}(b)),\cdots))$$
$$\times c(e_{2n'+1})\cdots
c(e_n)+ e^{-||b||^2}(\overline{\chi}<2n-2n')_b.\eqno(3.10)$$ As in
Section 2.3, using $\Omega^{\top}$ and $\Omega^{\bot}$ instead of
$A^{\top}$ and $A^{\bot}$ and multiplying the constant
$(-\sqrt{-1})^{[\frac{n}{2}]+1}2^{[\frac{n}{2}]}$, then we get
\begin{eqnarray*}
& & {\rm lim}_{t\rightarrow 0}\int_{\nu_{\xi}(\varepsilon)}{\rm
Tr}[\widetilde{\tau}^*(x)t^{\frac{k}{2}}f^0c(df^1)\cdots
c(df^k)\overline{P_t}(x)]dx\\
 &=&{\rm lim}_{t\rightarrow
0}\int_{{\bf R
}^{n-2n'}}(-\sqrt{-1})^{[\frac{n}{2}]+1}2^{[\frac{n}{2}]}
t^{\frac{n-2n'}{2}}t^{\frac{k}{2}}\\
& ~ & \times f^0(\xi)c(df^1)(\xi)\cdots c(df^k)(\xi)e^{-||b||^2}
\frac{(\sqrt{-1})^{m_1+1}}{(4\pi t)^{n/2}}\\
& ~ &\times
P(t;(\cdots,(2(-1)^l(\sum_{\alpha=1}^{n'}u^{2l}_{\alpha}+\sum_{\beta=1}^{m_q}v_{\beta}^{2l})),\cdots);
(\cdots,((-1)^l4t\sum_{\alpha=1}^{m_q}v_{\alpha}^{2l}(b^2_{2\alpha-1}+b^2_{2\alpha})),\cdots))db.
\end{eqnarray*}
Similar to the computation in Section 2.3, we obtain:

 \indent{\bf Theorem 3.6}~~{\it For $f^0,\cdots,f^k\in
 C^{\infty}_{\tau}(M)$ and $k$ even,}\\
$~~~{\rm lim}_{t\rightarrow 0}{\bf {\rm ch}}_k(\sqrt
{t}D)(f^0,\cdots ,f^k)$
 $$=\frac{1}{k!(2\pi\sqrt{-1})^{\frac{k}{2}}}
\sum_{q=1}^r \frac{(\sqrt{-1})^{m_1-m_q}}{2}\int_{F_q}f^0\wedge
df^1\wedge\cdots\wedge df^k\wedge \widehat{A}(TF_q)[{\rm
ch}\triangle(N_q)]^{-1}.\eqno(3.11)$$ {\it where $f^j$ is
considered as $f^j|_{F_q}$ for $0\leq j\leq k.$ }

\indent {\bf Remark:}~Since the computing of the Chern-Connes
character does not require the condition $f^j\in
 C^{\infty}_{\tau}(M)$, so (3.11) is correct for any $f^j\in
 C^{\infty}(M)$. When $k=0$ and $f^0=1$, we get the theorem 2.6.

{\footnotesize
 \indent {\bf Acknowledgements.}
 The author is indebted to
Professor Weiping Zhang for his guidance and very helpful
discussions. He also thanks Professor Huitao Feng for his generous
help and referees for their careful reading and helpful comments.}
\end{document}